\documentclass[a4paper,11pt]{article}
\usepackage[utf8]{inputenc}
\usepackage{a4wide}
\usepackage{latexsym,amsfonts,amsmath,amssymb,mathrsfs,url,amsthm}
\usepackage{mathtools}
\usepackage{color,graphicx}
\usepackage{subcaption}
\usepackage{caption}
\usepackage{lipsum}
\usepackage{xcolor}
\usepackage{hyperref}
\usepackage{verbatim}
\usepackage[normalem]{ulem}
\providecommand{\keywords}[1]
{
  \noindent \small	
  \textbf{Keywords:} #1
}

\providecommand{\amscode}[1]
{
  \noindent \small	
  \textbf{AMS subject classifications:} #1
}
\newtheorem{theorem}{Theorem}[section]

\newtheorem{remark}[theorem]{Remark}

\newcommand{\bsx}{\boldsymbol{x}}

\newcommand{\NN}{\mathbb{N}}

\newcommand{\ZZ}{\mathbb{Z}}

\DeclareMathOperator{\Var}{Var}
\DeclareMathOperator{\Unif}{Unif}
\DeclareMathOperator{\median}{median}

\title{The median trick does not help for fully nested scrambling\thanks{The work of T.G.\ is supported by JSPS KAKENHI Grant Number 23K03210. The work of K.~S.\ is supported by JSPS KAKENHI Grant Number 24K06857.}}
\author{Takashi Goda\thanks{Graduate School of Engineering, The University of Tokyo, 7-3-1 Hongo, Bunkyo-ku, Tokyo 113-8656, Japan (\url{goda@frcer.t.u-tokyo.ac.jp})}, Kosuke Suzuki\thanks{Faculty of Science, Yamagata University, 1-4-12 Kojirakawa-machi, Yamagata, 990-8560, Japan ({\tt kosuke-suzuki@sci.kj.yamagata-u.ac.jp})}}
\date{\today}

\begin{document}

\maketitle

\begin{abstract}
In randomized quasi-Monte Carlo methods for numerical integration, average estimators based on digital nets with fully nested and linear scrambling are known to exhibit the same variance.
In this note, we show that this equivalence does not extend to the median estimators.
Specifically, while the median estimator with linear scrambling can achieve faster convergence for smooth integrands, the median estimator with fully nested scrambling does not exhibit this advantage.
\end{abstract}


\amscode{65C05; 65D30; 65D32}

\section{Introduction}

In this paper, we consider the numerical integration problem over the unit cube. We approximate the integral
\[
I(f) := \int_{[0,1]^d} f(\bsx) \, d\bsx
\]
by averaging the function values at sampling nodes,
\[
Q_N(f) := \dfrac{1}{N} \sum_{i=0}^{N-1} f(\bsx_i).
\]
Among various approaches, randomized quasi-Monte Carlo (RQMC) is one of the most widely used methods for high-dimensional numerical integration.
RQMC combines the Monte Carlo (MC) method, which uses random sampling, with the quasi-Monte Carlo (QMC) method, which uses deterministic low-discrepancy point sets.
For details on MC, QMC, and RQMC, see \cite{DKP22,DKS13,DP10,Lem09,LL02,Nie92,SJ94}.

Typically, RQMC uses $r$ independent realizations $Q_N^{(i)}$ for $i = 1, \dots, r$, where each set of sampling nodes $\{\bsx_0^{(i)}, \dots, \bsx_{N-1}^{(i)}\}$ is independently and randomly drawn from some class of low-discrepancy point sets. Their average,
\begin{equation}\label{eq:average-estimator}
A_N^{(r)} := \dfrac{1}{r}\sum_{i=1}^r Q_N^{(i)}(f),
\end{equation}
is used as an estimator for $I(f)$. Therefore, the variance of $Q_N$ is critical in analyzing this estimator.

Recently, instead of the average $A_N^{(r)}$, the median of the realizations
\begin{equation}\label{eq:median-estimator}
M_N^{(r)} := \median(Q_N^{(1)}, \dots, Q_N^{(r)})
\end{equation}
has been intensively studied.
Remarkably, it has been shown that, when a suitable combination of low-discrepancy point sets and randomization methods is selected, the median estimator is \textit{universal}.
That is, it can achieve nearly optimal convergence rates across a variety of function spaces without prior knowledge of the function's smoothness \cite{GK24,GL22,GSM24,Pan25,PO23,PO24}.

In this note, we focus on one of the most well-known RQMC methods: scrambled nets.
In particular, we consider two specific types of scrambling: fully nested scrambling and linear scrambling (also known as affine matrix scrambling).
As outlined in \cite[Section~6.12]{DKS13}, see also \cite{Owen03}, it is known that scrambled net estimators, based on either fully nested scrambling or linear scrambling, exhibit the same variance.
Moreover, linear scrambling has been shown to benefit from the median trick: it is universal in the sense mentioned above \cite{GSM24,Pan25,PO23,PO24}.
This might lead one to conjecture that the median estimator for fully nested scrambling would share the same benefit.
The purpose of this note is to demonstrate that this is not the case:
the median estimator for fully nested scrambling fails to improve convergence rates.

\section{Results}
In this note, we work with the one-dimensional case.
For a fixed integer base $b\ge 2$, let $a = \sum_{k=1}^{\infty} \alpha_k b^{-k} \in [0,1)$ with $\alpha_k \in \ZZ_b := \{0,1,\dots,b-1\}$.
Scrambling randomizes the digits $\alpha_k$ to obtain $x = \sum_{k=1}^{\infty} \xi_k b^{-k}$,
where the digits $\xi_k$ are determined via permutations or matrix transformations.
The main scrambling types are summarized below, following \cite{Owen03}.

\paragraph{Fully nested scrambling}
In base $b$ fully nested scrambling \cite{Owen95}, each digit $\xi_k$ depends on preceding digits:
\[
\xi_k = \pi_{\bullet \alpha_1 \alpha_2 \dots \alpha_{k-1}}(\alpha_k),
\]
where the permutations $\pi_{\bullet \alpha_1 \alpha_2 \dots \alpha_{k-1}}$
are independently and randomly selected from all possible $b!$ permutations on $\ZZ_b$.
When $b$ is prime, one can use linear permutations $\pi(\alpha) = h \alpha + g \mod b$ instead.

\paragraph{Linear scrambling}
Let $\vec{a} = (\alpha_1, \alpha_2, \dots)^{\top}$ and $\vec{x} = (\xi_1, \xi_2, \dots)^{\top}$ be the digit vectors.
Matrix scrambling in base $b$ applies a lower-triangular matrix $M$ such that $\vec{x} = M \vec{a} \mod b$.
The matrix $M$ may take various forms, such as:
\[
M = \begin{bmatrix}
h_{11} & 0      & 0      & 0      & \cdots \\
g_{21} & h_{22} & 0      & 0      & \cdots \\
g_{31} & g_{32} & h_{33} & 0      & \cdots \\
g_{41} & g_{42} & g_{43} & h_{44} & \cdots \\
\vdots & \vdots & \vdots & \vdots & \ddots \\
\end{bmatrix},
\text{ or }
\begin{bmatrix}
h_1 & 0   & 0   & 0   & \cdots \\
g_2 & h_1 & 0   & 0   & \cdots \\
g_3 & g_2 & h_1 & 0   & \cdots \\
g_4 & g_3 & g_2 & h_1 & \cdots \\
\vdots & \vdots & \vdots & \vdots & \ddots \\
\end{bmatrix}
\text{ or }
\begin{bmatrix}
h_1 & 0   & 0   & 0   & \cdots \\
h_1 & h_2 & 0   & 0   & \cdots \\
h_1 & h_2 & h_3 & 0   & \cdots \\
h_1 & h_2 & h_3 & h_4 & \cdots \\
\vdots & \vdots & \vdots & \vdots & \ddots \\
\end{bmatrix}.
\]
Here, $h$'s are in $\ZZ_b^*:=\{1,\ldots,b-1\}$, and $g$'s are in $\ZZ_b$, drawn uniformly and independently.
We refer them as Matou\v{s}ek's matrix scrambling \cite{Mat98},
Tezuka’s $i$-binomial scrambling \cite{Tez94},
and Owen's striped matrix scrambling \cite{Owen03}, respectively.

\paragraph{Variance of scrambling}
For $m\in \NN$, let $\{a_0,\ldots,a_{b^m-1}\}\subset [0,1)$ be a $(0,m,1)$-net in base $b$, that is, a point set for which exactly one point lies in each of the disjoint intervals $[i/b^m, (i+1)/b^m)$ for $i=0,\ldots,b^m-1$.
The randomized estimator $Q_N$ is given by applying scrambling to each point $a_i$.
It is known that the resulting randomized point set $\{x_0,\ldots,x_{b^m-1}\}$ forms a $(0,m,1)$-net in base $b$ with probability one, and that the estimator $Q_N$ is unbiased.
For differentiable functions $f$ with Lipschitz continuous derivative of order $\beta$ (for a weaker assumption, see \cite{yue1999variance}),
both scrambling types yield the same variance \cite{Owen1997b}:
\[
\Var(Q_N) = \frac{1}{12N^3} \int_{0}^1 (f'(x))^2 \,dx +  O(N^{-3-\beta}),\qquad \text{with $N=b^m$}.
\]
Moreover, for fully nested scrambling, the central limit theorem holds \cite{Loh03}:
\begin{equation}\label{eq:CLT-Q}
N^{3/2}(Q_N-I(f)) \xrightarrow{d} \mathcal{N}(0, \sigma^2(f))
\qquad \text{with } \sigma^2(f) = \dfrac{1}{12}\int_{0}^1 (f'(x))^2 \,dx.
\end{equation}

\paragraph{Median estimator of the linear scrambling}

The median estimator for linearly scrambled nets has been studied recently in \cite{GSM24,Pan25,PO23,PO24}.
In the one-dimensional setting considered here, the most relevant result is given in \cite{PO23}, where Pan and Owen showed that the median estimator automatically exploits the smoothness of the integrand---without requiring it to be specified---and achieves an improved convergence rate with high probability, specifically faster than $1/N^{3/2}$.
For infinitely smooth functions, it can even attain super-polynomial convergence, a convergence behavior originally shown in \cite{Suz17}.
Furthermore, for functions with finite smoothness, the result has been strengthened in \cite{Pan25}, where improved rates are obtained when the number of replicates $r$ increases logarithmically with $N$.

A heuristic explanation is that, for smoother functions, the error distribution of linearly scrambled nets becomes more heavy-tailed, while the variance remains the same as in \eqref{eq:CLT-Q}.
As a result, taking the median over multiple replicates significantly reduces the probability that the final estimate is influenced by ``outliers'' in the error distribution.

\paragraph{Median estimator of the fully nested scrambling}
In one dimension, fully nested scrambling in base $b$ of a $(0,m,1)$-net is equivalent to jittered sampling: the $i$-th point $X_{N,i}$ is uniformly sampled from $[i/N, (i+1)/N]$.
Thus, the estimator becomes
\begin{equation}\label{eq:full-estimator}
Q_N = \frac{1}{N} \sum_{i=0}^{N-1} f(X_{N,i}),
\qquad \text{where $X_{N,i} \overset{\mathrm{iid}}{\sim} \Unif{\left[\dfrac{i}{N}, \dfrac{i+1}{N}\right]}$ and $N=b^m$.}
\end{equation}
The main result of this note is the following asymptotic behavior of the median estimator, which is proven in Section~\ref{sec:proof}. 
The proof is a straightforward application of some basic theorems from order statistics.

\begin{theorem}\label{thm:asymptotic}
Let $f$ have Lipschitz continuous derivative of order $\beta$, and define
\[
I(f) := \int_0^1 f(x)\,dx, \qquad \sigma^2(f) := \dfrac{1}{12}\int_{0}^1 (f'(x))^2\,dx.
\]
Let $Q_N$ be the estimator of the fully nested scrambling given in \eqref{eq:full-estimator}
and $M_N^{(r)}$ the median estimator using $r$ samples of $Q_N$'s,
defined as in \eqref{eq:median-estimator}.
Then we have the following convergence in law:
\[
\sqrt{r}\,N^{3/2}(M_{N}^{(r)} - I(f))
\xrightarrow{d}
\mathcal{N}\left(0,\dfrac{\pi\,\sigma^{2}(f)}{2}\right)
\qquad (N \to \infty \text{ and then } r \to \infty).
\]
\end{theorem}

\begin{remark}
This result implies that, even for functions smoother than required in the theorem, the asymptotic error distribution remains normal.
In contrast, for linearly scrambled nets, the error distribution becomes increasingly heavy-tailed as the smoothness of the integrand increases.
As a consequence of this observation, the median estimator based on fully nested scrambling does not yield an improved convergence rate for smooth functions, unlike that based on linear scrambling. 
\end{remark}

In median RQMC, the number of realizations $r$ is typically fixed or grows slowly with respect to $N$.
For example, setting $r=15$ is justified by a rule of thumb \cite[Remark~4.1]{GSM24}.
It has been proven in \cite{GK24,Pan25} that a slowly increasing function, such as $r(N) = \log N$,
ensures that the median estimator retains an almost optimal randomized error bound across a wide class of function spaces.

The median estimator based on Mato\v{u}sek's linear scrambling has been shown to be universal.
For functions $f \in C^{\alpha}[0,1]$ with $\alpha\ge 1$,
the root-mean-squared error achieves a rate of $O(N^{-\alpha - 1/2 + \epsilon})$, for any $\epsilon > 0$, without requiring knowledge of the smoothness parameter $\alpha$. 

In contrast, for such reasonable (i.e., constant or slowly increasing) choices of $r$,
Theorem~\ref{thm:asymptotic} strongly suggests that the error of the median estimator based on the fully nested scrambling is typically of order $N^{-3/2}r(N)^{-1/2}$ even if $f$ is very smooth.
Thus, the fully nested scrambling cannot exploit the smoothness of integrands.

Fortunately, most existing implementations of scrambling methods are based on linear scrambling.
Nevertheless, the result presented here highlights the need for caution when combining scrambling techniques with the median trick.

\section{Numerial experiments}
Consider the following two functions:
\[ f_1(x)=x^{3/2},\]
and 
\[ f_2(x)=\exp(-x).\]
The first function is finitely smooth, whereas the second is infinitely smooth.
This implies that, with a proper choice of $r$ depending on $N$, the median estimator based on linear scrambling achieves an error decay of order $N^{-\alpha-1/2}$ for some $\alpha>1$ when applied to $f_1$, and of order $N^{-c\log N}$ for some $c>0$ for $f_2$.
This convergence behavior, however, is not expected for the median estimator based on fully nested scrambling.

\begin{figure}[t]
  \centering
  \begin{subfigure}[t]{0.45\textwidth}
    \centering
    \includegraphics[width=\linewidth]{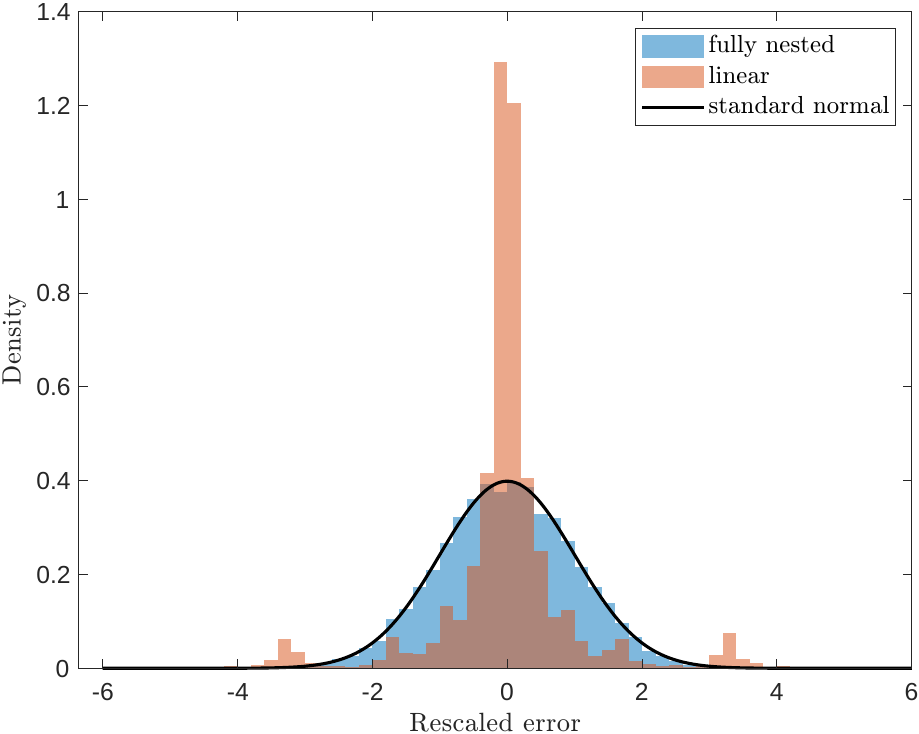}
    \caption{$f_1$, $N=16$}
  \end{subfigure}
  \hfill
  \begin{subfigure}[t]{0.45\textwidth}
    \centering
    \includegraphics[width=\linewidth]{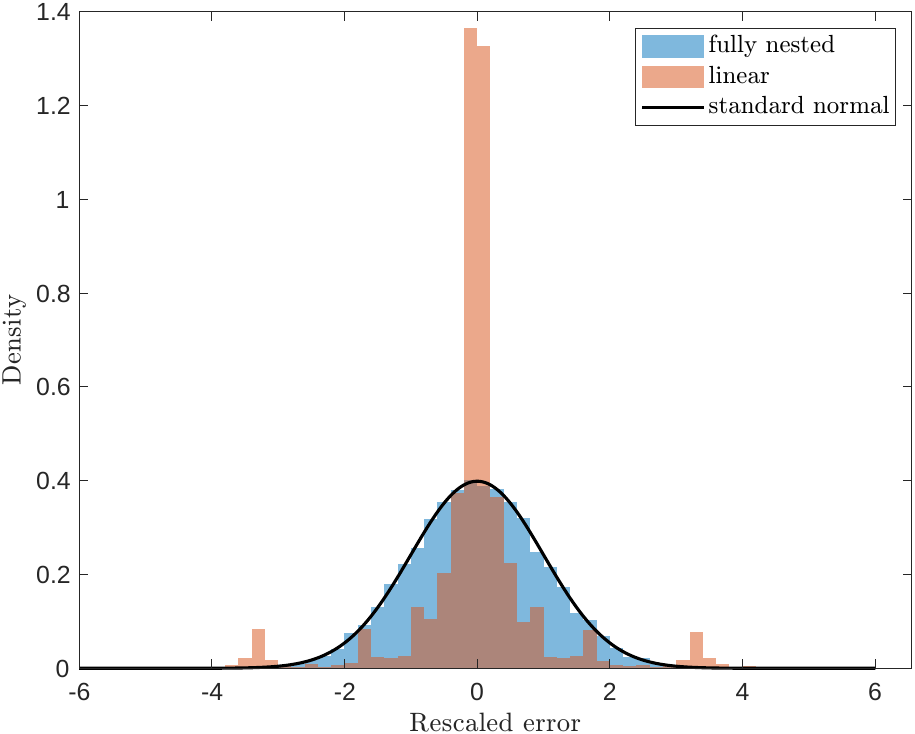}
    \caption{$f_2$, $N=16$}
  \end{subfigure}

  \vspace{0.5cm}

  \begin{subfigure}[t]{0.45\textwidth}
    \centering
    \includegraphics[width=\linewidth]{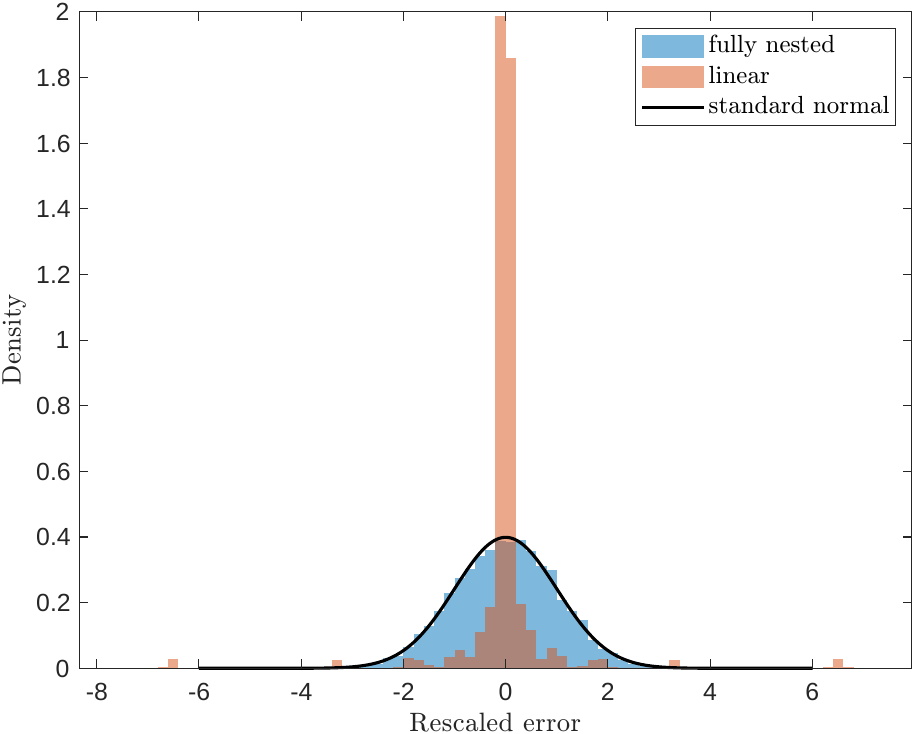}
    \caption{$f_1$, $N=64$}
  \end{subfigure}
  \hfill
  \begin{subfigure}[t]{0.45\textwidth}
    \centering
    \includegraphics[width=\linewidth]{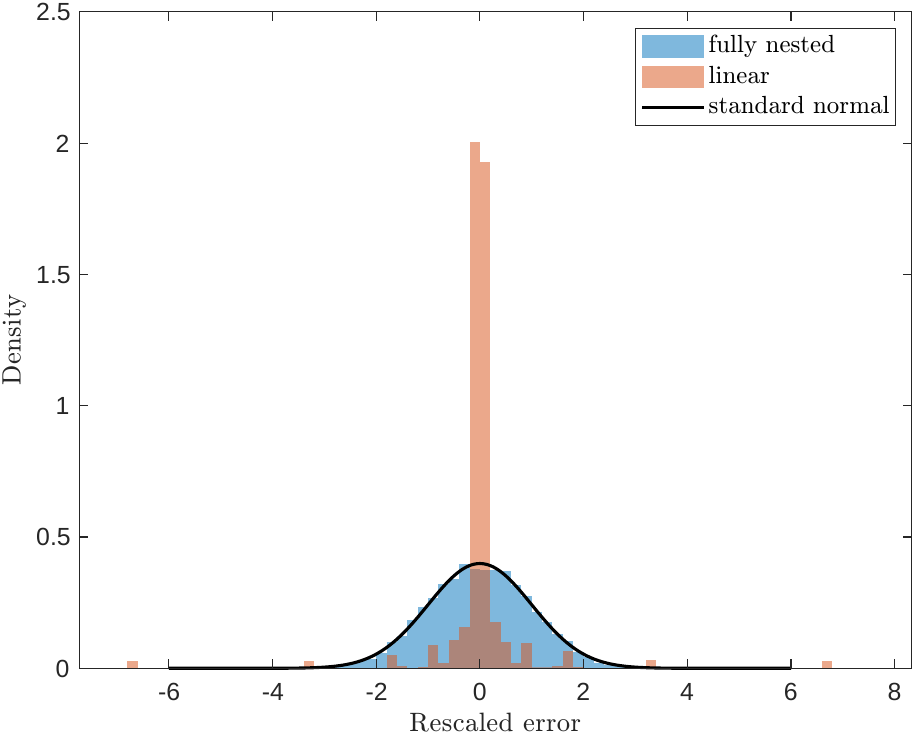}
    \caption{$f_2$, $N=64$}
  \end{subfigure}

  \caption{Histograms of rescaled errors for different test functions and values of $N$ with $r=1$.}
  \label{fig:histograms}
\end{figure}

First, to observe the difference in the error distributions between fully nested scrambling and linear scrambling, we compute the estimates $Q_N$ using each scrambling method $10^4$ times and compare their histograms in Figure~\ref{fig:histograms}. In each subplot, the horizontal axis shows the rescaled error $N^{3/2}(Q_N(f_i) - I(f_i))/\sigma(f_i)$, and the vertical axis shows the empirical density. We also plot the probability density function of the standard normal distribution as a reference. As expected from the literature \cite{Loh03,Owen03}, the rescaled error for fully nested scrambling closely follows the standard normal distribution. Although the error distributions differ significantly between the two scrambling methods, we checked that the empirical variances are nearly identical and match the theoretical variance, which is equal to $1$ as expected from \eqref{eq:CLT-Q}.

\begin{figure}[t]
  \centering
  \begin{subfigure}[t]{0.45\textwidth}
    \centering
    \includegraphics[width=\linewidth]{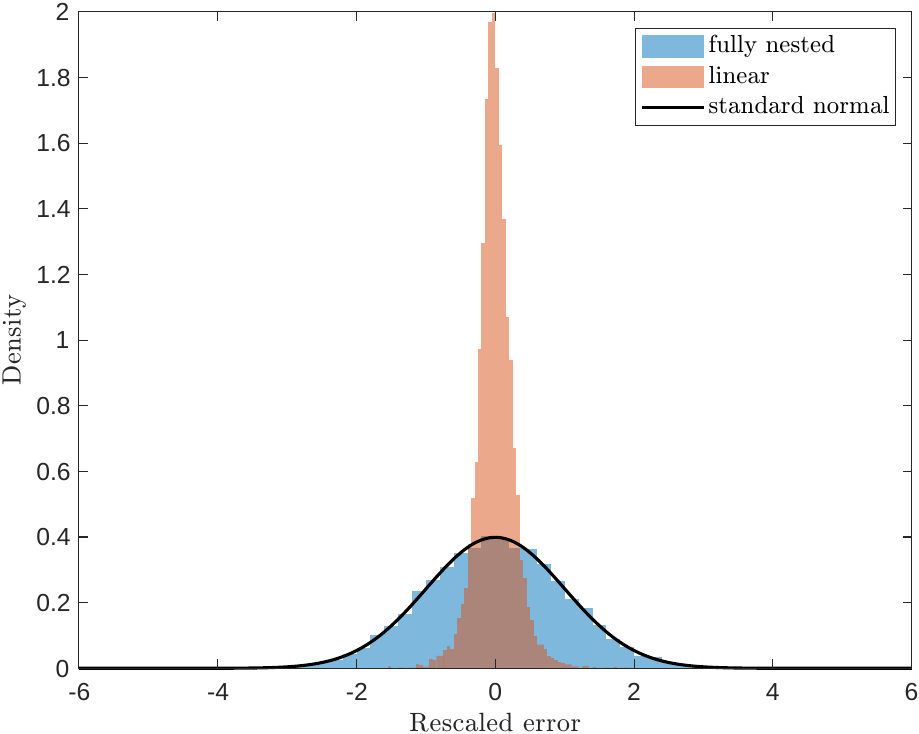}
    \caption{$f_1$, $N=16$}
  \end{subfigure}
  \hfill
  \begin{subfigure}[t]{0.45\textwidth}
    \centering
    \includegraphics[width=\linewidth]{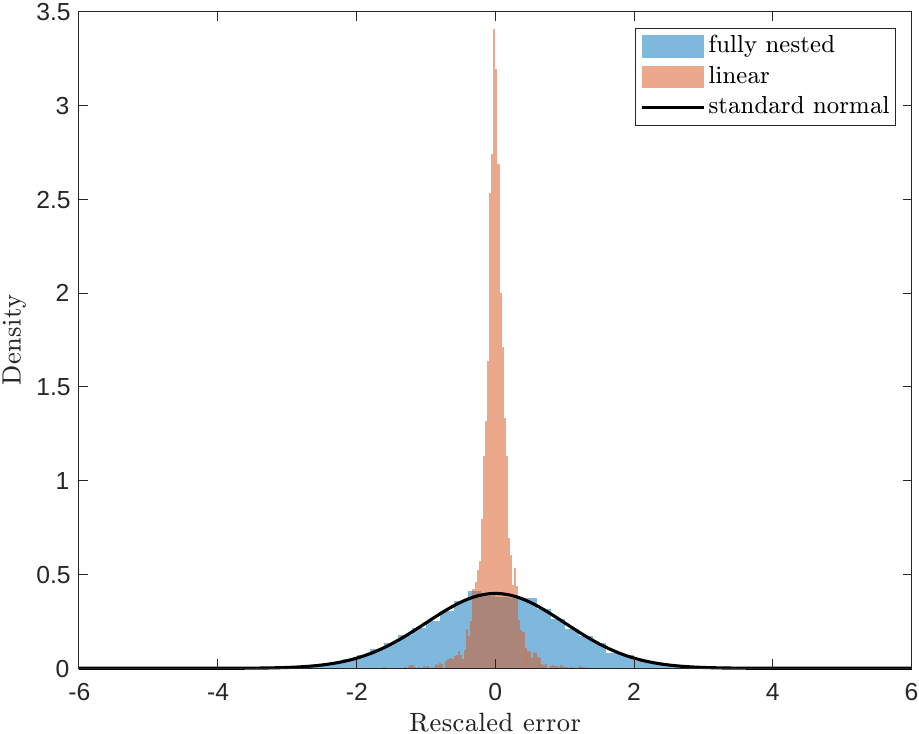}
    \caption{$f_2$, $N=16$}
  \end{subfigure}

  \vspace{0.5cm}

  \begin{subfigure}[t]{0.45\textwidth}
    \centering
    \includegraphics[width=\linewidth]{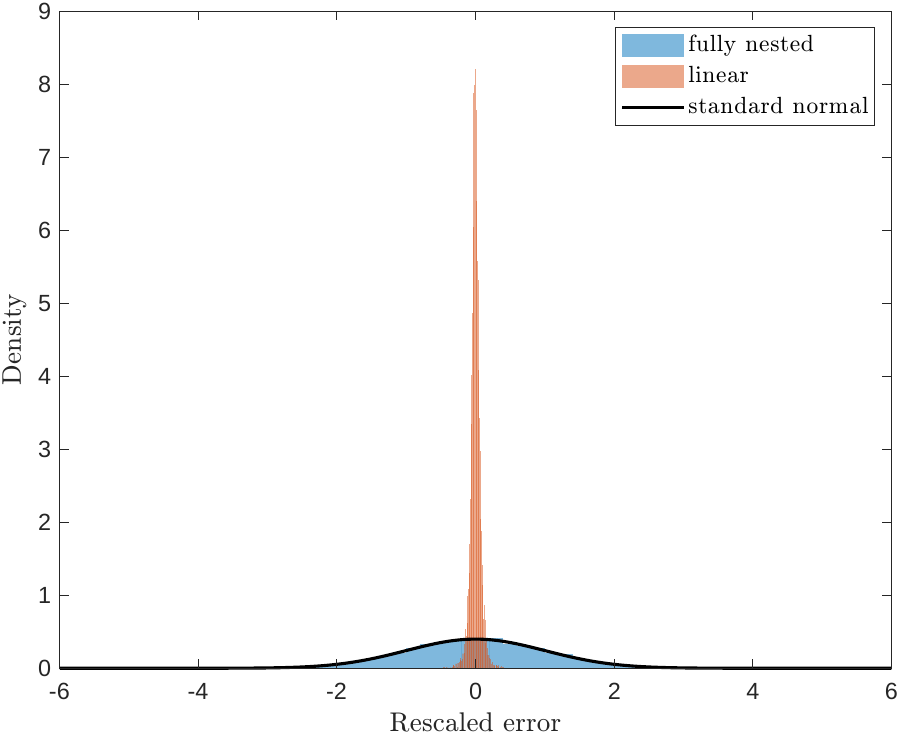}
    \caption{$f_1$, $N=64$}
  \end{subfigure}
  \hfill
  \begin{subfigure}[t]{0.45\textwidth}
    \centering
    \includegraphics[width=\linewidth]{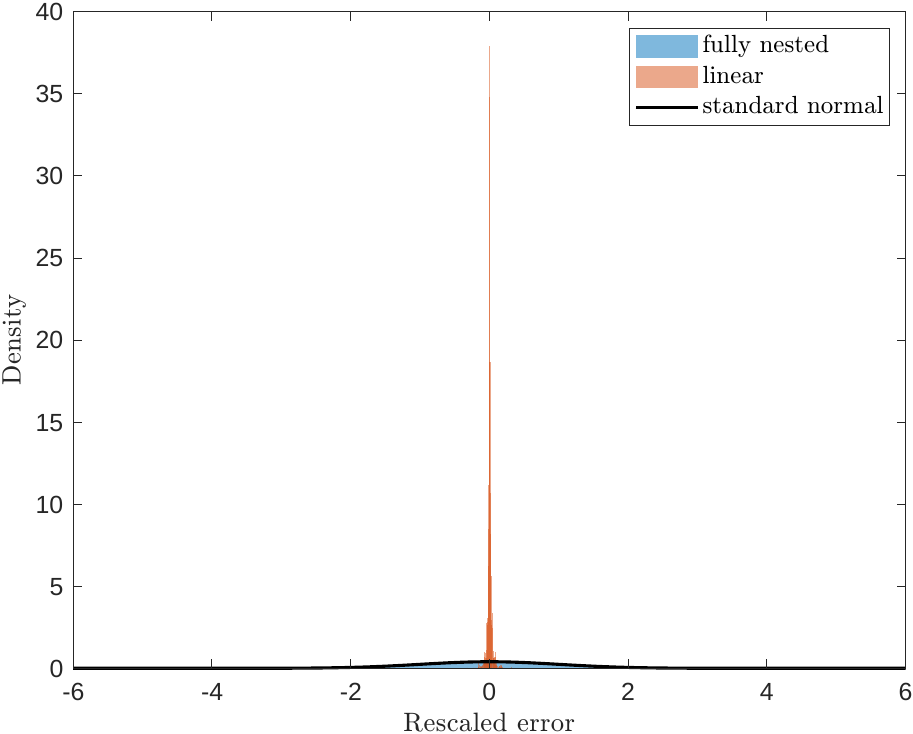}
    \caption{$f_2$, $N=64$}
  \end{subfigure}

  \caption{Histograms of rescaled errors for different test functions and values of $N$ with $r=15$.}
  \label{fig:histograms_r15}
\end{figure}

Next, we compare the error distributions of the median estimator $M_N^{(15)}$, based on $r=15$ independent estimates, for fully nested scrambling and linear scrambling. As in the above experiment, we compute $M_N^{(15)}$ $10^4$ times and present the histograms in Figure~\ref{fig:histograms_r15}. In each subplot, the horizontal axis shows the rescaled error $\sqrt{2r/\pi}N^{3/2}(M_N^{(r)} - I(f_i))/\sigma(f_i)$, and the vertical axis shows the empirical density. We also plot the probability density function of the standard normal distribution as a reference. As expected from our main result (Theorem~\ref{thm:asymptotic}), the rescaled error for fully nested scrambling again follows the standard normal distribution. In contrast, the rescaled error for linear scrambling is quite different; it is more concentrated around zero, resulting in smaller errors. Moreover, the variance of the rescaled error decreases as $N$ increases, which agrees well with the fact that linear scrambling can exploit the smoothness of integrands. This effect is especially pronounced for $f_2$, which has infinite smoothness, compared to $f_1$.

\begin{figure}[t]
  \centering
  \begin{subfigure}[t]{0.45\textwidth}
    \centering
    \includegraphics[width=\linewidth]{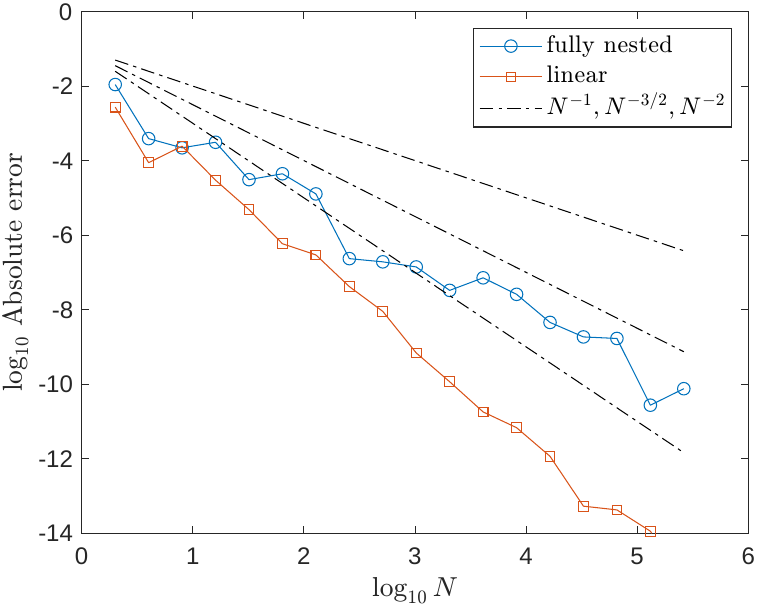}
    \caption{$f_1$}
  \end{subfigure}
  \hfill
  \begin{subfigure}[t]{0.45\textwidth}
    \centering
    \includegraphics[width=\linewidth]{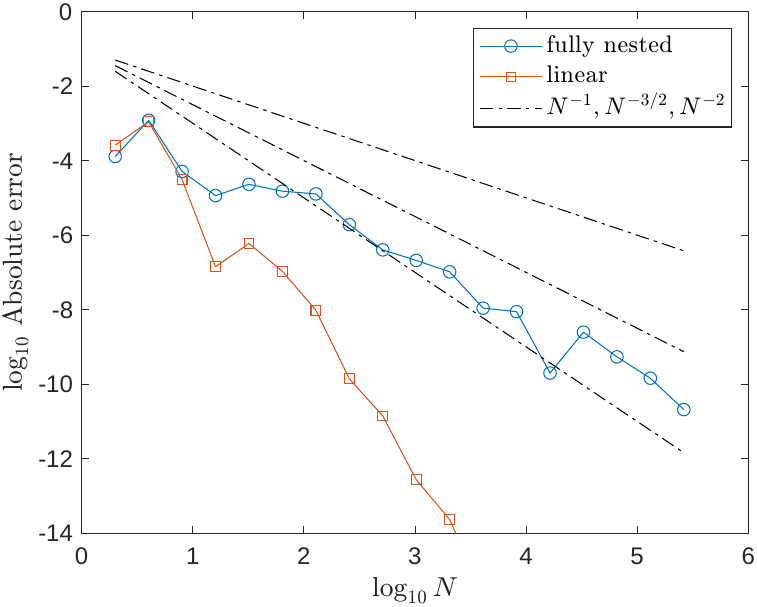}
    \caption{$f_2$}
  \end{subfigure}
  \caption{Convergence behavior of the median estimators $M_N^{(r)}$ for both fully nested scrambling and linear scrambling.}
  \label{fig:convergence}
\end{figure}

Figure~\ref{fig:convergence} shows the convergence behavior of the median estimators $M_N^{(r)}$ for both fully nested scrambling and linear scrambling. To mitigate fluctuations in the median estimators, we set $r=1001$. In each subplot, the horizontal axis indicates $\log_{10} N$, where $N=2^m$ is the number of sampling nodes, and the vertical axis shows $\log_{10}$ of the absolute error of each median estimator. The three dashed lines correspond to the rates $\mathcal{O}(N^{-1})$, $\mathcal{O}(N^{-3/2})$, and $\mathcal{O}(N^{-2})$ for reference.
For fully nested scrambling, the convergence order is consistent with the rate $\mathcal{O}(N^{-3/2})$, as expected from Theorem~\ref{thm:asymptotic}. This gives empirical support that the median estimator based on fully nested scrambling does not achieve an improved convergence rate. In contrast, for linear scrambling, the convergence is significantly faster, particularly for the smoother integrand $f_2$. This observation is consistent with known theoretical results in the literature, which state that linear scrambling can exploit the smoothness of integrands, leading to faster error decay.

\section{Proof of the main theorem}\label{sec:proof}
We end this note with a proof of Theorem~\ref{thm:asymptotic}.

We first consider the asymptotic distribution of the sample median as $N \to \infty$ for fixed $r$.
Let $Q_N$ be the estimator corresponding to fully nested scrambling.
Let $r$ be an odd integer, and form $r$ independent replicates $\{Q_{N}^{(j)}\}_{j=1}^{r}$,
each distributed identically to $Q_N$.

Define the sample median
$M_{N}^{(r)}$ as in \eqref{eq:median-estimator} and let
\begin{equation}
\widetilde{M}_{N}^{(r)}     
:= N^{3/2}(M_{N}^{(r)} - I(f))
= \median\{S_{N}^{(1)},\dots,S_{N}^{(r)}\}.
\end{equation}
where
\[
S_{N}^{(j)} := N^{3/2}(Q_N^{(j)} - I(f)).
\]
From \eqref{eq:CLT-Q}, we have 
\[
S_{N}^{(j)}\xrightarrow{d}\mathcal{N}(0,\sigma^{2}) \qquad (N \to \infty).
\]
Thus, it follows from the continuous-mapping theorem that
\[
\median\{S_{N}^{(1)},\dots,S_{N}^{(r)}\} 
\xrightarrow{d}
\median\{Z_{1},\dots,Z_{r}\} \qquad (N \to \infty),
\]
where $Z_i \sim \mathcal{N}(0, \sigma^2)$ are i.i.d.\ for $i=1,\dots,r$.
Since $r = 2k+1$ is odd, the density function of $\median\{Z_1, \dots, Z_r\}$ is known from the theory of order statistics \cite[Chapter~2]{ABN08}, given by
\[
f_{\mathrm{med}}^{(r)}(x)
:= \frac{(2k-1)!}{(k-1)!\,(k-1)!}
[\Phi_{\sigma}(x)]^{\,k-1}
[1 - \Phi_{\sigma}(x)]^{\,k-1}
\,\varphi_{\sigma}(x),
\]
where $\Phi_{\sigma}$ and $\varphi_{\sigma}$ are the CDF and PDF of $\mathcal{N}(0,\sigma^2)$, respectively.
Consequently, $\widetilde{M}_{N}^{(r)}$ converges in distribution to a random variable with density $f_{\mathrm{med}}^{(r)}(x)$.

Now consider the case where $r \to \infty$ as well.
It is known that the sample median is asymptotically unbiased and normal, and that, for $r \to \infty$, we have
\[
\sqrt{r} \cdot \median\{Z_{1},\dots,Z_{r}\}
\xrightarrow{d}
\mathcal{N}\left(0,\dfrac{\pi\,\sigma^{2}}{2}\right),
\]
see for instance \cite[Chapter~8]{ABN08}.
Thus, taking $N \to \infty$ first and then $r \to \infty$, we obtain
\[
\sqrt{r}\,\widetilde{M}_{N}^{(r)}
\xrightarrow{d}
\mathcal{N}\left(0,\dfrac{\pi\,\sigma^{2}}{2}\right),
\]
which establishes Theorem~\ref{thm:asymptotic}.

%%%%%%%%%%%%%%%%%%%%%%%%%%%%%%
%%%%%%%%%%%%%%%%%%%%%%%%%%%%%%
\bibliographystyle{siam}
\bibliography{ref.bib}

\end{document}